# Estimating the Evolution of Solution Norms in Vector Delay Nonlinear Systems: Stability and Boundedness

Mark A. Pinsky

***Abstract*—** Existing methods rarely capture the temporal evolution of solution norms in vector nonlinear DDEs with variable delays and coefficients, often leading to overly conservative boundedness and stability criteria. We develop a framework that constructs scalar counterparts of vector DDEs whose solutions upper-bound the evolution of the original solution norms when the corresponding history functions are matched. This reduction enables boundedness and stability assessment of vector DDEs through the dynamics of their scalar counterparts, using straightforward simulations or simplified analytical reasoning. New boundedness and stability criteria and the estimates of the radii of balls containing history functions that yield bounded or stable solutions for the original vector systems were derived and validated through representative simulations.

## I. Introduction

Assessing the boundedness and stability of vector nonlinear systems with variable delays and coefficients is a long-standing problem with broad applications in science and engineering. Extensive studies reported in monographs [4], [7]-[9], [11]-[15], [17]-[20], [24] and surveys [1], [5], [10], [16] have proposed various stability and some boundedness criteria for specific classes of delay systems.

Recent advances primarily stem from extensions of the Lyapunov methodology for time-delay systems introduced by Krasovskii [19] and Razumikhin [29]. These developments have been further complemented by the use of Halanay's inequality [7], [8], [12], [22], the Bohl–Perron theorem [2], [3], [9], and several related more specialized approaches. For linear time-invariant systems, such methods have led to computationally efficient LMI techniques [8], [18], [30], while stability of time-varying DDEs has been studied in [3], [4],[6], [21], [33].

However, applying such methods to nonlinear DDEs with variable delays and coefficients remains restrictive, producing overly conservative criteria that rarely provide acceptable estimates of stability or boundedness regions [23], [31]. These limitations are especially pronounced in nonhomogeneous vector nonlinear systems of this type, where ISS methods address certain aspects of the problems [5], but many challenges remain unresolved.

Current stability and boundedness paradigm focuses on the system's asymptotic behavior but disregard the transient dynamics, which may be crucial in practical applications.

Building on our previous work on nonlinear ODEs [25]–[28], this study develops a new framework for analyzing the temporal evolution of solution norms in vector nonlinear systems with variable delays and coefficients. The method constructs scalar auxiliary counterparts for the original vector systems. Under appropriate matching of history functions, solutions of these scalar equations upper-bound the evolution of solution norms of vector DDEs.

This reduction, enabling estimation of the evolution of solution norms of vector DDEs through simplified analytical reasoning and efficient simulations of their scalar counterparts, constitutes the main contribution of this paper. Consequently, we derive new boundedness and stability criteria and estimate for radii of the balls embedded in boundedness or stability regions. Linearization of the scalar auxiliary equations under Lipschitz-like condition leads to simplified, albeit more conservative criteria. The proposed approach is validated by representative simulations that confirm its accuracy and practical relevance.

The paper is organized as follows. Section II introduces notation, definitions and system formulation; Section III derives the scalar auxiliary equation and the associated boundedness–stability criteria; Section IV develops the linearized scalar auxiliary equations together with the corresponding boundedness–stability criteria. Section V reports simulation results, and Section VI concludes and remarks on future research directions.

## II. Notation and Governing Equations

**2.1. Notation.** Recall that the symbols $\mathbb{R}$, $\mathbb{R}_{\geq 0}$, $\mathbb{R}_+$ and $\mathbb{R}^n$ represent the sets of real, non-negative, and positive real numbers and real $n$-dimensional vectors, $\mathbb{N}$ is a set of real positive integers, $\mathbb{R}^{n \times n}$ is a set of $n \times n$ matrices, and $I \in \mathbb{R}^{n \times n}$ is an identity matrix. Next, $C([a,b];\mathbb{R}^n)$, $C([a,b];\mathbb{R}_+)$, and $C([a,b];\mathbb{R}_{\geq 0})$ are the spaces of continuous real functions, $\zeta:[a,b] \to \mathbb{R}^n$, $\zeta:[a,b] \to \mathbb{R}_+$ and $\zeta:[a,b] \to \mathbb{R}_{\geq 0}$, respectively, with the supremum norm $\|\zeta\| := \sup_{t \in [a,b]} |\zeta(t)|$, where $|\cdot|$ stands for the Euclidean norm of a vector or the induced norm of a matrix, $a < b \in \mathbb{R}$ and $b$ can be infinite. In turn, $\dot{x}(t) := D^+ x(t)$, where $D^+$ is the upper-right-hand derivative in $t$. Also note that $|x|_\infty = \sup_{i=1,\ldots,n}(|x_i|)$ and $|x|_1 = \sum_{i=1}^{n} |x_i|$, $\forall x \in \mathbb{R}^n$.

**2.3. Underlined Equation and Definitions.**
The subsequent analysis concerns the vector nonlinear equation with variable delays and coefficients:

*Mark A. Pinsky with the Department of Mathematics & Statistics, University of Nevada.Reno, Reno NV 89557 USA (e-mail: pinsky@unr.edu)*

$$\dot{x} = A(t)x + f(t, x(t), x(t-h_1(t)), \ldots, x(t-h_m(t))) + F(t),$$
$$\forall t \geq t_0; \; x(t, \varphi) = \varphi(t), \; \forall t \in [t_0 - \bar{h}, t_0] \quad (2.1)$$

where $x(t, \varphi) \in \Re \subset \mathbb{R}^n, 0 \in \Re, f(t, \chi_1, \ldots, \chi_{m+1}) \in \mathbb{R}^n$, $\forall t \geq t_0, \forall \chi_i \in \mathbb{R}^n, i = 1, \ldots, m+1$ is a continuous function in all variables, $f(t, 0) = 0$, the initial function $\varphi \in C([t_0 - \bar{h}, t_0]; \mathbb{R}^n), \|\varphi\| := \sup |\varphi(t)|, t \in [t_0 - \bar{h}, t_0]$, scalar functions $h_i \in C([t_0, \infty); \mathbb{R}_+), \sup h_i(t) \leq \bar{h} < \infty$, matrix $A \in C([t_0, \infty); \mathbb{R}^{n \times n})$ is continuously differentiable, $F(t) = F_0 e(t), \; e \in C([t_0, \infty); \mathbb{R}^n), \|e\| = \sup |e(t)| = 1, t \geq t_0$, $F_0 \in \mathbb{R}_{\geq 0}$.

We assume that the initial problem (2.1) admits a unique solution $\forall t \geq t_0$, and $\forall \|\varphi\| \leq \bar{\varphi} > 0$. Under the above assumptions, the right-hand side of (2.1) is continuous in all variables, ensuring that a solution to this equation, $x(t, t_0, \varphi)$, is continuous, continuously differentiable in $t, \forall t \geq t_0$, and bounded on any finite interval. Thus, the ensuing analysis estimates for both the transient evolution of the solution norms of (2.1) and their asymptotic behavior as $t \to \infty$.

Henceforth, we adopt the abridged notation: $x(t, \varphi) := x(t, t_0, \varphi), \forall t \geq t_0$ and introduce the homogeneous counterpart of (2.1):
$$\dot{x} = A(t)x + f(t, x(t), x(t-h_1(t)), \ldots, x(t-h_m(t))),$$
$$\forall t \geq t_0, \; x(t, \varphi) = \varphi(t), \; \forall t \in [t_0 - \bar{h}, t_0] \quad (2.1)$$

which admits admit a unique solution $\forall t \geq t_0$. We also consider the associated linear equation:
$\dot{x} = A(t)x, \forall t \geq t_0, x(t_0, \varphi) = \varphi(t_0) \in \mathbb{R}^n$ whose solution can be expressed as: $x(t, \varphi(t_0)) = W(t, t_0)\varphi(t_0)$, where $W(t, t_0) = w(t)w^{-1}(t_0)$ is the transition (Cauchy) matrix, and $w(t)$ is a fundamental solution matrix of the linear equation. The latter formular defines the variation of parameters in the next section.

Then, we assume that $|w(t_0)| = 1$, a condition that is always satisfiable.

Next, under the above-stated assumptions, we recall the standard definitions of boundedness and stability of solutions to equations (2.1) and (2.2), which will be used in the sequel; see, for instance, [18].

**Definition 1**. The trivial solution of equation (2.2) is called: (1.1) Stable for the given value of $t_0$ if $\forall \varepsilon \in \mathbb{R}_+, \exists \delta_1(t_0, \varepsilon) \in \mathbb{R}_+$ such that $|x(t, t_0, \varphi)| < \varepsilon, \forall t \geq t_0$, $\forall \|\varphi\| < \delta_1(t_0, \varepsilon)$. Otherwise, it is unstable. (1.2.) Uniformly stable if in the above definition $\delta_1(t_0, \varepsilon) = \delta_2(\varepsilon)$. (1.3) Asymptotically stable if it is stable for given value of $t_0$ and $\exists \delta_3(t_0) \in \mathbb{R}_+$ such that $\lim_{t \to \infty} |x(t, t_0, \varphi)| = 0, \quad \forall \|\varphi\| < \delta_3(t_0)$.
(1.4) Uniformly asymptotically stable if it is uniformly stable and in the previous definition $\delta_3(t_0) = \delta_4 = const$.

**Definition 2**. The solution to equation (2.1) is called: (2.1) Locally bounded for the fixed value of $t_0$ if $\exists \delta_6(t_0) \in \mathbb{R}_+$ and $\exists \varepsilon_*(\delta_6) \in \mathbb{R}_+$ such that a condition $\forall \|\varphi\| < \delta_6(t_0)$ implies that $|x(t, t_0, \varphi)| < \varepsilon_*, \forall t \geq t_0$. (2.2) Uniformly bounded if $\delta_6(t_0) = \delta_7 = const$.

Furthermore, let $R_i \equiv \sup \delta_i$ be the superior value of $\delta_i$ for which the conditions of subdefinitions (1.3) or (1.4) of Definition 1 or conditions of subdefinitions (2.1) or (2.2) of Definition 2 hold. Next, let $B_{R_i} := \{\varphi \in C([t_0 - \bar{h}, t_0], \mathbb{R}^n) : \|\varphi\| \leq R_i\}, i = 3, \ldots, 6$ denote centered at the origin balls with radius $R_i$. These balls contain history functions generating bounded or stable solutions of equations (2.1) and (2.2), respectively.

Next, we recall a definition of robust stability, termed stability under persistent perturbations by Krasovskii [31], pp.161-164. Consider a perturbed form of equation (2.2):
$$\dot{x} = A(t)x + f(t, x, x(t-h_1(t)), \ldots, x(t-h_m(t))) +$$
$$R(t, x, x(t-h_1^*(t)), \ldots, x(t-h_m^*(t))), \forall t \geq t_0, \quad (2.2)$$
$$x(t) = \varphi(t), \; \forall t \in [t_0 - \bar{h}, t_0]$$

where in addition to the conditions specified above for the components of (2.2) included in (2.3), we also assume that: (I) function $f(t, \chi_1, \ldots, \chi_{m+1})$ is Lipschitz continuous in $\chi_i \in \mathbb{R}^n, i = 1, \ldots, m+1$, (II) $f(t, 0) = 0$, (III) $R(t, \chi_1, \ldots, \chi_{m+1}) \in \mathbb{R}^n, \forall t \geq t_0, \forall \chi_i \in \mathbb{R}^n, i = 1, \ldots, m+1$ is a continuous function in all variables, though $R(t, 0)$ may be nonzero, and (IV) scalar functions $h_i^* \in C([t_0, \infty); \mathbb{R}_+)$.

**Definition 3**. The trivial solution to equation (2.2) is called robustly stable if $\forall \varepsilon \in \mathbb{R}_+, \exists \Delta_i(\varepsilon) \in \mathbb{R}_+, i = 1, 2, 3$ such that $|x(t, \varphi)|_\infty \leq \varepsilon, \forall t \geq t_0$ if
$|R(t, \chi_1, \ldots, \chi_{m+1})|_\infty < \Delta_1(\varepsilon), \forall t \geq t_0, \forall |\chi_i|_\infty < \varepsilon, i = 1, \ldots, m+1$,
$\|\varphi\| < \Delta_2(\varepsilon), |h_i(t) - h_i^*(t)| < \Delta_3(\varepsilon), \forall t \geq t_0, i = 1, \ldots, m$, and $x(t, \varphi)$ is a solution of equation (2.3).

This definition characterizes the robust stability of (2.2) through the boundedness of the solutions of the perturbed equation (2.3).

**Statement** (Krasovskii [19], p.162). Under conditions (I)–(IV), the trivial solution of (2.2) is robustly stable if it is uniformly asymptotically stable.

### III. Scalar Delay Auxiliary Equations

A scalar auxiliary DDE with a history function is derived below to upper-bound the solution norms of vector DDE (2.1).

Both equations are cast in integral form via the use of variation of parameters. Applying standard norms and extended Lipschitz inequalities yields to a scalar integral equation for (2.1) that is matched then to the integral form of scalar auxiliary DDE.

First, we write (2.1) in the integral form as:

$$x(t,\varphi) = w(t)w^{-1}(t_0)\varphi(t_0) + w(t)\int_{t_0}^{t} w^{-1}(\tau)Q(\tau)d\tau, \quad \forall t \geq t_0,$$

$$x(t,\varphi) = \varphi(t), \quad \forall t \in [t_0 - \bar{h}, t_0]$$

where $Q = f(t, x(t), x(t-h_1(t)), \ldots, x(t-h_m(t))) + F(t)$.

Consequently, we get that:

$$|x(t,\varphi)| = \left| w(t)w^{-1}(t_0)\varphi(t_0) + w(t)\int_{t_0}^{t} w^{-1}(\tau)Q(\tau)d\tau \right|, \quad (3.1)$$

$$\forall t \geq t_0, |x(t,\varphi)| = |\varphi(t)|, \quad \forall t \in [t_0 - \bar{h}, t_0]$$

Then, applying standard norm inequalities yields:

$$|X_1(t,\varphi)| = |w(t)||w^{-1}(t_0)\varphi(t_0)| + |w(t)|\int_{t_0}^{t} |w^{-1}(\tau)||Q(\tau)|d\tau, \quad (3.1)$$

$$\forall t \geq t_0, |X_1(t,\varphi)| = |\varphi(t)|, \quad \forall t \in [t_0 - \bar{h}, t_0]$$

where $X_1(t,\varphi) \in \mathbb{R}^n$ is a continuous vector-function in both variables. Comparing (3.1) and (3.2) implies that $|x(t,\varphi)| \leq |X_1(t,\varphi)|, \forall t \geq t_0$.

To render a more tractable form of (3.2), we introduce a nonlinear extension of the Lipschitz condition:

$$|f(t,\chi_1,\ldots,\chi_{m+1})| \leq L(t,|\chi_1|,\ldots,|\chi_{m+1}|),$$

$$\forall \chi = [\chi_1,\ldots,\chi_{m+1}]^T \in \Omega \in \mathbb{R}^{n(m+1)}, \chi_i \in \mathbb{R}^n \quad (3.2)$$

where $\Omega$ is a compact subset of $\mathbb{R}^{n(m+1)}$ containing zero, $f \in C([t_0,\infty) \times \mathbb{R}^{n(m+1)}; \mathbb{R}^n)$, $f(t,0) = 0$, and $L(t,0) = 0$, $L \in C([t_0,\infty) \times \mathbb{R}_{\geq 0}^{m+1}; \mathbb{R}_{\geq 0})$ is a scalar continuous function.

Appendix B, illustrates how to define such scalar functions in a closed form if $f(t,\chi_1,\ldots,\chi_{m+1})$ is a polynomial or power series in $\chi_1,\ldots,\chi_{m+1}$ with a bounded in $\Omega$ error term. In the former case, $\Omega \equiv \mathbb{R}^{n(m+1)}$ and the same condition also holds if power series error term is bounded in $\mathbb{R}^{n(m+1)}$. Moreover, $L$ becomes linear in $|\chi_i|$ if $f$ is linear in the corresponding variables. Furthermore, Section V defines $L$ for a planar nonlinear system.

Note that the remainder of this paper assumes for simplicity that $\Omega \equiv \mathbb{R}^{n(m+1)}$.

Next, substituting (3.3) into (3.2) yields the following equation:

$$X_2(t,\tilde{\varphi}) = |w(t)||w^{-1}(t_0)\varphi(t_0)| + |w(t)|\int_{t_0}^{t} |w^{-1}(\tau)|q(\tau)d\tau, \quad (3.3)$$

$$\forall t \geq t_0; X_2(t,\tilde{\varphi}) = |\varphi(t)|, \quad \forall t \in [t_0 - \bar{h}, t_0]$$

where $X_2(t,\tilde{\varphi}) \in \mathbb{R}_{\geq 0}$ is a continuous scalar function of both variables and

$q = L(t, X_2(t), X_2(t-h_1(t)), \ldots, X_2(t-h_m(t))) + |F(t)|$.

Hence, equation (3.4) implies that $|x(t,\varphi)| \leq |X_1(t,\varphi)| \leq X_2(t,\tilde{\varphi}), \forall t \geq t_0$.

Next, we introduce an auxiliary scalar nonlinear DDE:

$$\dot{y} = p(t)y + c(t)\big(L(t,y(t),y(t-h_1(t)),\ldots,y(t-h_m(t))) + |F(t)|\big),$$

$$\forall t \geq t_0, \quad y(t,\phi) = \phi(t), \quad \forall t \in [t_0 - \bar{h}, t_0] \quad (3.5)$$

where, $y(t,\phi) \in \mathbb{R}_{\geq 0}$, $\phi \in C([t_0 - \bar{h}, t_0]; \mathbb{R}_{\geq 0})$, and functions $p:[t_0,\infty) \to \mathbb{R}$, $c:[t_0,\infty) \to [1,\infty]$, and $\phi$ are determined below. Subsequently, we write (3.5) in the integral form using variation of parameters:

$$y(t) = e^{d(t)}\left(\phi(t_0) + \int_{t_0}^{t} e^{-d(\tau)} c(\tau) q(\tau) d\tau\right), \quad \forall t \geq t_0 \quad (3.4)$$

$$y(t,\phi) = \phi(t), \quad \forall t \in [t_0 - \bar{h}, t_0]$$

where $d(t) = \int_{t_0}^{t} p(s)ds$.

To determine $p(t)$ and $c(t)$, we match the right-hand sides of (3.6) and (3.4). Matching the first terms: $|w(t)||w^{-1}(t_0)\varphi(t_0)|$ and $e^{d(t)}\phi(t_0)$, returns:

$$|w(t)| = \exp\left(\int_{t_0}^{t} p(s)ds\right) \quad (3.5)$$

$$\phi(t_0) = |w^{-1}(t_0)\varphi(t_0)| \quad (3.6)$$

To match the second terms on the right side of (3.6) and (3.4), we first multiply and divide the latter function by $|w(t)|$ and use that, due to (3.7), $e^{-d(t)} = 1/|w(t)|$, which yields:

$$\int_{t_0}^{t} e^{-d(\tau)} c(\tau) q(\tau) d\tau = \int_{t_0}^{t} \big(|w(\tau)||w^{-1}(\tau)|\big) q(\tau)/|w(\tau)| d\tau$$

The last relation shows that $c(t) = |w(t)||w^{-1}(t)|$ is the running condition number of $w(t)$. Since, $0 < |w^{-1}(t)| < \infty, \forall t \geq t_0$, it follows that $c(t) < \infty, \forall t \geq t_0$. Furthermore, the continuous differentiability of $A(t)$ ensures the same property for $c(t)$ and continuity of $p(t)$.

Lastly, (3.7) yields that:

$$p(t) = d(\ln|w(t)|)/dt \quad (3.7)$$

Matching the initial functions in (3.4) and (3.6) yields that $\phi(t) = |\varphi(t)|, \forall t \in [t_0 - \bar{h}, t_0]$, which combined with (3.8) implies that: $|w^{-1}(t_0)\varphi(t_0)| = |\varphi(t_0)|$. This condition, together with previous requirement, $|w(t_0)| = 1$, is satisfied upon setting $w(t_0) = I$.

Thus, (3.5), with the defined above coefficients and supplemented by the condition:

$y(t,\phi) = \phi(t) = |\varphi(t)|$, $\forall t \in [t_0 - \bar{h}, t_0]$, specifies the initial problem for a scalar auxiliary DDE whose solutions upper-bound the norms of the corresponding solutions of the original vector equation (2.1). This leads to the following statement.

**Theorem 1.** Assume that $f \in C([t_0, \infty) \times \mathbb{R}^{n(m+1)}; \mathbb{R}^n)$, $f(t, 0) = 0$, matrix $A \in C([t_0, \infty); \mathbb{R}^{n \times n})$ is continuously differentiable, $\varphi \in C([t_0 - \bar{h}, t_0]; \mathbb{R}^n)$, $h_i \in C([t_0, \infty); \mathbb{R}_+)$, $\sup h_i(t) \leq \bar{h} < \infty$, $F(t) = F_0 e(t)$, $e \in C([t_0, \infty); \mathbb{R}^n)$, $\|e\| = 1$, $F_0 \in \mathbb{R}_{\geq 0}$, $L \in C([t_0, \infty) \times \mathbb{R}_{\geq 0}^{m+1}; \mathbb{R}_{\geq 0})$, $L(t, 0) = 0$, inequality (3.3) holds with $\Omega \equiv \mathbb{R}^{n(m+1)}$, and $w(t_0) = I$. Subsequently, we assume that equations (2.1) and (3.5) assume unique solutions $\forall \|\varphi\| \leq \bar{\varphi} > 0$, $\forall t \geq t_0$. Then,

$$|x(t,\varphi)| \leq y(t,\phi), \forall t \geq t_0 \quad (3.8)$$
$$\phi(t) = |\varphi(t)|, \forall t \in [t_0 - \bar{h}, t_0]$$

where $x(t,\varphi)$ and $y(t,\phi)$ are the solutions to equations (2.1) and (3.5), respectively.

**Proof.** In fact, it has been shown previously that $|x(t,\varphi)| \leq X_2(t, \tilde{\varphi})$, $\forall t \geq t_0$. Assignments on the functions $p(t)$, $c(t)$ and $\phi(t)$ made above yield: $X_2(t, \tilde{\varphi}) \equiv y(t, \phi)$, $\forall t \geq t_0$ if $\tilde{\varphi}(t) = \phi(t) = |\varphi(t)|$, $\forall t \in [t_0 - \bar{h}, t_0]$ □

Therefore, a boundedness or stability condition imposed on the scalar equation (3.5) or on its homogeneous counterpart obtained for $F_0 = 0$ ensures that the same property carries over to vector equations (2.1) and (2.2), respectively.

Assume next that Definitions 1 and 2 are applied to solutions of scalar equation (3.5) or to those of its homogeneous analog with $F_0 = 0$. Then, define $r_i \equiv \sup \delta_i$, $i = 3,...,6$ for which the conditions of subdefinitions (1.3) or (1.4) in Definition 1, when applied to (3.5) with $F_0 = 0$, or conditions of subdefinitions (2.1) or (2.2) in Definition 2, when applied to (3.5) with $F_0 > 0$, hold. In turn, let $B_{r_i} := \{\varphi \in C([t_0 - \bar{h}, t_0], \mathbb{R}^n) : \|\varphi\| \leq r_i\}$, $i = 3,...,6$ be the balls with radiuses $r_i$ that are centered at the origin.

This leads to the following statements.

**Theorem 2.** Assume that the conditions of Theorem 1 are satisfied and that the trivial solution to equation (3.5) with $F_0 = 0$ is either stable, uniformly stable, asymptotically stable, or uniformly asymptotically stable. Then, the trivial solution to equation (2.2) exhibits the same type of stability.

Furthermore, under the conditions (1.3) and (1.4) of Definition 1, $B_{r_i} \subseteq B_{R_i}$, $i = 3, 4$.

**Proof.** Both above statements directly follow from the application of inequality (3.10) to the corresponding solutions of equations (2.2) and (3.5) with $F_0 = 0$ □

**Theorem 3.** Assume that the conditions of Theorem 1 are met and that the solutions to equation (3.5) are bounded for the set values of $t_0$ if $\|\varphi\| \leq r_5$ or uniformly bounded if $\|\varphi\| \leq r_6$. Then, the solutions to equation (2.1) with the matched history functions are also bounded for the same values of $t_0$ or are uniformly bounded. Furthermore, under prior conditions, $B_{r_i} \subseteq B_{R_i}$, $i = 5, 6$.

**Proof.** The proof of this statement follows directly from the application of inequality (3.10) to the solutions of equation (2.1) and its scalar counterpart (3.5) □

The key task of estimating the values of $r_i$, $i = 3,...,6$ is abridged by the following statement.

**Lemma 1.** Assume that $y_i(t, \phi_i)$, $\forall t \geq t_0$, $i = 1, 2$ are the solutions to (3.5) and $\phi_1(t) \geq \phi_2(t)$, $\forall t \in [t_0 - \bar{h}, t_0]$. Then, $y_1(t, \phi_1) \geq y_2(t, \phi_2)$, $\forall t \geq t_0$.

**Proof.** The proof of this statement follows from the assumption of uniqueness of solutions to equation (3.5), which rules out the intersection of the solution curves to a scalar delay equation in $t \times y$ - plane □

Thus, a solution to equations (3.5) and of its homogeneous counterpart with a constant history function, i.e, $y = y(t, \phi)$, $\forall t \geq t_0$, $\phi(t) \equiv q \in \mathbb{R}$, $\forall t \in [t_0 - \bar{h}, t_0]$, is monotonic in $q$, $\forall t \geq t_0$. This further streamlines the simulations of the values of $r_i = \sup \|\varphi\|$, $i = 3,...,6$ for which the corresponding conditions of Definitions 1 and 2 hold.

Next, we apply our approach to ensure the robust stability of the trivial solution of the vector delay equation (2.2) via conditions imposed on its scalar counterpart. Consequently, we write the scalar counterpart for (2.3) as follows:

$$\dot{z} = p(t)z + c(t)\begin{pmatrix} L(t, z(t), z(t-h_1(t)),...,z(t-h_m(t))) + \\ L_R(t, z(t), z(t-h_1^*(t)),...,z(t-h_m^*(t))) \end{pmatrix}, \quad (3.9)$$

$\forall t \geq t_0$; $z(t, \phi) = |\varphi(t)|$, $\forall t \in [t_0 - \bar{h}, t_0]$

where $z(t, \phi) \in \mathbb{R}_{\geq 0}$ and $L_R \in C([t_0, \infty) \times \mathbb{R}_{\geq 0}^{m+1}; \mathbb{R}_{\geq 0})$.

Then, the application of Definition 3 to equations (3.5) with $F_0 \equiv 0$ and (3.11) leads to the following statement.

**Theorem 4.** Assume that the conditions of Theorem 1 and conditions (I)–(IV) hold for equation (3.11), and that the trivial solution of the homogeneous counterpart of equation (3.5) is uniformly asymptotically stable. Then, this trivial solution and, consequently, that of equation (2.2), are robustly stable.

**Proof.** In fact, under the above conditions, the trivial solution of homogeneous counterpart of (3.5) is robustly stable which implies that $z(t, |\varphi|) < \varepsilon$, $\forall t \geq t_0$ (where $z(t, |\varphi|)$ is a solution to (3.11)) if $|L_R(t, \xi_1,...,\xi_{m+1})|_\infty < \Delta_1(\varepsilon)$, $\forall t \geq t_0$, $\forall \xi_i < \varepsilon$, $i = 1,..., m+1$, $\|\varphi\| < \Delta_2(\varepsilon)$, and $|h_i(t) - h_i^*(t)| < \Delta_3(\varepsilon)$, $\forall t \geq t_0$, see Definition 3. Then, (3.11) yields that $|x(t, \varphi)| \leq z(t, |\varphi|) < \varepsilon$, $\forall t \geq t_0$, where $x(t, \varphi)$ is a solution to (2.3), since (3.11) is the scalar counterpart to (2.3) □

**Corollary.** Assume that the conditions of Theorem 1 hold, and $\sup_{\forall t \geq t_0} p(t) = \hat{p} < 0$, $\sup_{\forall t \geq t_0} c(t) = \hat{c} < \infty$, and $\sup_{\forall t \geq t_0} L(t, y) = \hat{L}(y) < \infty$. Next, we assume that $\hat{p} y + \hat{c} \hat{L}(y) < 0$, $\forall y \in (0, y_+)$, $y_+ > 0$ and conditions (I)–(IV) hold for equation (3.11). Then, the trivial solution of the corresponding vector equation (2.2) is robustly stable.

**Proof.** Under the given assumptions, the trivial solution to the corresponding scalar equation (3.5) with $F_0 \equiv 0$ is uniformly asymptotically stable. By Theorem 4, this implies that equation (2.2) is robustly stable □

Thus, the proposed approach upper-bound evolution of solution norms of vector nonlinear DDEs with variable delay and coefficients, thereby streamlining boundedness and stability analysis of such equations through rigorously justified simulations of the associated scalar counterparts. These simulations also estimate the radii of history-function balls leading to bounded or asymptotically stable solutions of the original system. Moreover, the approach simplifies the application of analytical tools to these problems.

## IV. Linearized Scalar Delay Auxiliary Equations

Under a Lipschitz-type condition, this section derives simplified yet more conservative criteria for the boundedness and stability of scalar nonlinear delay equations and their associated vector counterparts.

Let us assume that,

$$L(t, \zeta_1, ..., \zeta_{m+1}) \leq \sum_{i=1}^{m+1} \mu_i(t, \tilde{\zeta}) \zeta_i, \quad \forall |\zeta| \leq \tilde{\zeta} > 0, \quad \forall t \geq t_0 \quad (4.1)$$

where $\zeta_i \in \mathbb{R}_{\geq 0}$, $\zeta = [\zeta_1, ..., \zeta_{m+1}]^T$ and $\mu_i \in C([t_0, \infty) \times \mathbb{R}_{\geq 0}; \mathbb{R}_{\geq 0})$. Then, the application of (4.1) to (3.5) yields a scalar linear equation:

$$\dot{u} = P(t) u + c(t) \left( \sum_{i=2}^{m+1} \mu_i(t, \tilde{\zeta}) u(t - h_i(t)) + F_0 |e(t)| \right)$$
$$u(t, \phi) = |\varphi(t)|, \quad \forall t \in [t_0 - \bar{h}, t_0] \quad (4.2)$$

where $P(t) = p(t) + c(t) \mu_1(t, \tilde{\zeta})$. Next, we set in (4.1) that,

$$\zeta_1 = u(t, \phi), \zeta_i = u(t - h_i(t), \phi), i = 2, ..., m+1 \quad (4.3)$$

where $u(t, \phi)$ is a solution to either equation (4.2) or its homogeneous counterpart obtained for $F_0 = 0$. Consequently, we formulate the following stability criterion.

**Theorem 5.** Assume that $\mu_i \in C([t_0, \infty) \times \mathbb{R}_{\geq 0}; \mathbb{R}_{\geq 0})$, $i = 1, ..., m+1$, $p \in C([t_0, \infty); \mathbb{R})$, $c \in C([t_0, \infty); \mathbb{R}_+)$, $\varphi \in C([t_0, \infty); \mathbb{R}^n)$, and $h_i(t) \in C([t_0, \infty); \mathbb{R}_{>0})$, satisfy condition (A.1), see Appendix A, and that equation (4.2) admits a unique solution $\forall \tilde{\zeta} < \zeta_* > 0$, $\forall t \geq t_0$, and $\forall \|\varphi\| \leq \hat{\varphi} > 0$. Next, let the conditions of subdefinitions (1.i), i=1,...,4 of Definition 1 hold for the trivial solution of equation (4.2) with $F_0 \equiv 0$ and $\forall \tilde{\zeta} \leq \tilde{\zeta}_{i,s} \leq \zeta_*$, where $\tilde{\zeta}_{i,s} = \sup \tilde{\zeta}$ for which the conditions of subdefinition (1.i) hold. Assume further that the conditions of Theorem 1 are met. Then,

$$|x(t, \varphi)| \leq y(t, \phi) \leq u(t, \phi), \forall t \geq t_0 \quad (4.4)$$

where $x(t, \varphi)$, $y(t, \phi)$ and $u(t, \phi)$ are the solutions to equations (2.2) and homogeneous counterparts of (3.5) and (4.2), respectively. Moreover, the trivial solutions to equations (4.2) and (3.5) with $F_0 \equiv 0$, as well as of equation (2.2), exhibits the same type of stability as defined by subdefinition (1.i) of Definition 1.

**Proof.** Inequality (4.1) holds under the least restrictive conditions of subdefinition (1.1), which implies that $u(t, \phi) < \varepsilon$, $\forall t \geq t_0$ if $\|\phi\|$ is sufficiently small, where $u(t, \phi)$ solves (4.2) with $F_0 \equiv 0$. This ensures that (4.1) is satisfied for sufficiently small $\varepsilon$ and $\forall \tilde{\zeta} \leq \tilde{\zeta}_{1,s}$. Then, (4.1) ensures applicability of Lemma 2 (see Appendix A) to homogeneous counterparts of (3.5) and of (4.2), which yields that $y(t, \phi) \leq u(t, \phi)$, $\forall t \geq t_0$, conforming the stability of both the trivial solutions of (3.5) with $F_0 \equiv 0$ and of (2.2) for the set value $t_0$.

Furthermore, the same conclusions follow under the more conservative conditions of subdefinitions (1,i), $i = 2, 3, 4$ □

In turn, the solution of (4.2) can be represented as (see, e.g., [18], p.139):

$$u(t, \phi) = u_h(t, \phi) + F_0 u_{nh}(t) \quad (4.5)$$

where $u_{nh}(t)$ is a solution to (4.2) with $\varphi(t) \equiv 0$ and $F_0 = 1$, and $u_h(t, \phi)$ solves (4.2) with $F_0 \equiv 0$. Furthermore, $u_{nh}(t) = \int_{t_0}^{t} C(t, s) |e(s)| ds$, $\forall t \geq t_0$, where $C(t, s)$ is the Cauchy function of (4.2) defined as a solution to homogeneous counterpart of (4.2) with the history function: $C(t, s) = \begin{cases} 0, & t_0 - \bar{h} \leq t < s \\ 1, & t = s \end{cases}$. This yields the following boundedness criterion.

**Theorem 6.** Assume that the conditions of Theorem 5 hold and equation (4.2) admits a unique solution under these conditions. Further, assume that $\sup |u_{nh}(t)| < \infty$, $t \geq t_0$, and both $F_0$ and $\|\varphi\|$ are sufficiently small. Then,

$$|x(t, \varphi)| \leq y(t, \phi) \leq u(t, \phi), \forall t \geq t_0 \quad (4.6)$$

where $x(t, \varphi), y(t, \phi)$, and $u(t, \phi)$ are the solutions to equations (2.1), (3.5), and (4.2), respectively, and $\phi(t) = |\varphi(t)|$, $\forall t \in [t_0 - \bar{h}, t_0]$.

**Proof.** Under the stated assumptions, both terms on the right-hand side of equation (4.5) can be made arbitrarily small by the appropriate choice of $\|\varphi\|$ and $F_0$, which implies that (4.1) holds. Hence, due to (4.1), Lemma 2 is applicable to equations (3.5) and (4.2), thereby ensuring that inequality (4.6) holds □

Thus, relations (4.5) and (4.6) imply that $|x(t,\varphi)| \leq u_h(t,|\varphi(t)|) + F_0 u_{nh}(t)$, where, under Theorem 6, $|u_{nh}(t)| < \infty, \forall t \geq t_0$ and $u_h(t,\phi)$ is stable or asymptotically stable solution of the homogeneous counterpart of (4.2), depending upon the adopted criterion. The condition $\sup|u_{nh}(t)| < \infty, \forall t \geq t_0$ eliminates input–system resonance in scalar linear equation (4.2) and ensures input-to-state stability of vector nonlinear DDE (2.2) under all stability conditions specified in Definition 1. This approach avoids Lyapunov-like arguments and leads to a new ISS condition (see [5] for further discussion).

## V. Simulations

Let us assume that (2.1) takes the following form,

$$\dot{x} = (A_0(t) + A_1(t))x + \rho A_1(t)x(t-h) + f(t,x(t-h)) + F(t), \forall t \geq t_0 = 0, x(t,x_0) \in \mathbb{R}^2, x(t,x_0) = x_0 = const, \forall t \in [-h,0] \quad (5.1)$$

where, $f = b[0, x_2^3(t-h)]^T$, $F = [0, F_2(t)]^T$, $b, h, \rho \in \mathbb{R}_+$,

$$A_0 = diag(\lambda(t), \lambda(t)), A_1 = \begin{pmatrix} 0 & 1 \\ -\omega(t) & -\alpha_1 \end{pmatrix}.$$

Next, we calculate the functions $p(t)$ and $c(t)$ for the subsystem $\dot{x} = A_0(t)x$ (see Section 3). Clearly, in this case, $w(t) = diag(\eta, \eta)$, where $\eta = \exp\left(\int_{t_0}^{t} \lambda(s)ds\right)$, which implies that $p(t) = \lambda(t)$ and $c(t) \equiv 1$.

Thus, for (5.1), the scalar auxiliary equation (3.5) becomes:

$$\dot{y} = (\lambda(t) + |A_1(t)|)y + \rho|A_1(t)|y(t-h) + by^3(t-h) + |F(t)|, \forall t \geq 0, y(t) \in \mathbb{R}_{\geq 0}; y(t) = |x_0|, \forall t \in [-h,0] \quad (5.2)$$

In turn, its autonomous scalar counterpart ends up as:

$$\dot{\hat{y}} = (\hat{\lambda} + \hat{A}_1)\hat{y} + \rho\hat{A}_1\hat{y}(t-h) + b\hat{y}^3(t-h) + F_0,$$

$\forall t \geq 0, \hat{y}(t) \in \mathbb{R}_{\geq 0}; \hat{y}(t) = |x_0|, \forall t \in [-h,0]$

(5.3) where $\sup_{\forall t \geq t_0}\lambda(t) = \hat{\lambda} < 0$, $\sup_{\forall t \geq t_0}|A_1(t)| = \|A_1\| = \hat{A}_1$. Thus, Lemma 2 implies that $y(t) \leq \hat{y}(t), \forall t \geq 0$.

Recall that the stability of the linearized equation (5.3) is well understood in the space of its parameters, and its asymptotic stability guarantees the asymptotic stability of the trivial solution of the homogeneous counterpart of equation (5.3), provided that $b$ is sufficiently small [18]. Moreover, the right–hand sides of equations (5.2) and (5.3) depend only on the scalar function $|A_1(t)|$ and the scalar parameter $|A_1(t)|$, respectively, rather than on all four entries of the matrix $|A_1(t)|$. This substantially simplifies the sensitivity analysis of the original equation (5.1) with respect to parameter variations.

Simulations of equations (5.1)-(5.3) assume that $\lambda = \lambda_0 + \lambda_+(t)$, where $\lambda_0 = -3$ and either (a)

$\lambda_+(t) = q\sin(dt), q = 0.1, d = 5$ or (b)

$\lambda_+(t) = q\exp(-dt), q = 1, d = 1$. Next, we assume that $\omega(t) = 1 + \omega_1(t), \omega_1(t) = 0.1(\sin(t) + \sin(\pi t))$ and $e(t) = [0, \sin(10t)]^T$. All simulations were caried out using the MATLAB solver DDE23.

Figure 1 compares the temporal evolution of the solution norms for equations (5.1)-(5.3), simulated with matched history functions and shown as solid, dashed, and dash-dotted lines, respectively. All parameters are indicated on the figure.

In all plots, inequality $|x(t,x_0)| \leq y(t,|x_0|) \leq \hat{y}(t,|x_0|)$ holds over simulated time intervals and across broad parameter ranges, supporting our theoretical inferences. Solutions of the scalar nonautonomous equation (5.2) closely approximate the solution norms of (5.1), whereas simulations of the autonomous equation (5.3) yield coarser estimates.

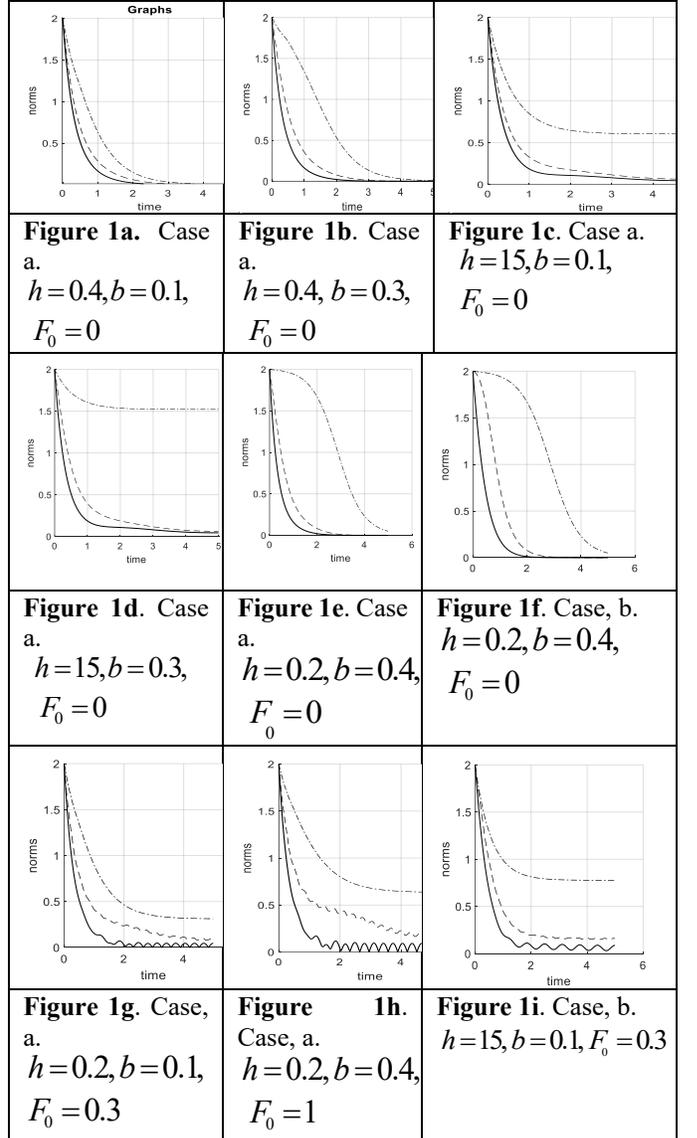

| Figure 1a. Case a. $h=0.4, b=0.1$, $F_0=0$ | Figure 1b. Case a. $h=0.4, b=0.3$, $F_0=0$ | Figure 1c. Case a. $h=15, b=0.1$, $F_0=0$ |
| Figure 1d. Case a. $h=15, b=0.3$, $F_0=0$ | Figure 1e. Case a. $h=0.2, b=0.4$, $F_0=0$ | Figure 1f. Case, b. $h=0.2, b=0.4$, $F_0=0$ |
| Figure 1g. Case, a. $h=0.2, b=0.1$, $F_0=0.3$ | Figure 1h. Case, a. $h=0.2, b=0.4$, $F_0=1$ | Figure 1i. Case, b. $h=15, b=0.1, F_0=0.3$ |

**Figure 1**. Evolution of solution norms for the vector and scalar DDEs.

To simulate the boundaries of the trapping and stability regions, defined as sets of initial vectors leading to bounded or stable solutions of equation (5.1) or its homogeneous counterpart, the initial vector was represented in polar coordinates. The angle was discretized with a step $\pi/100$, and for each angle, the radius increased until a sharp rise in $|x(t,x_0)|$, indicated that the boundary has been crossed, was observed. The same approach was applied to scalar equations (5.2) and (5.3) to estimate disk radii within true trapping or stability regions. The computation time for these estimates is nearly independent of the number of coupled equations and further decreases due to monotonic growth of $y(t,|x_0|)$ and $\hat{y}(t,|x_0|)$ in $|x_0|$, $\forall t \geq t_0$. In contrast, boundary estimation for vector equation (2.1) scales with $O(m^n)$, where $m$ is the number of discretization points per phase-space variable and $n$ is the number of such variables.

Figure 2 simulates the boundaries of the trapping ($F_0 > 0$) and stability ($F_0 = 0$) regions for equations (5.1)-(5.3) shown in solid, dashed, and dashed dotted lines for parameter sets indicated on the plots. Dashed and dash-dotted curves nearly coincide on all figures. The plots use $\ln r$ instead of $r$ to enhance visual clarity.

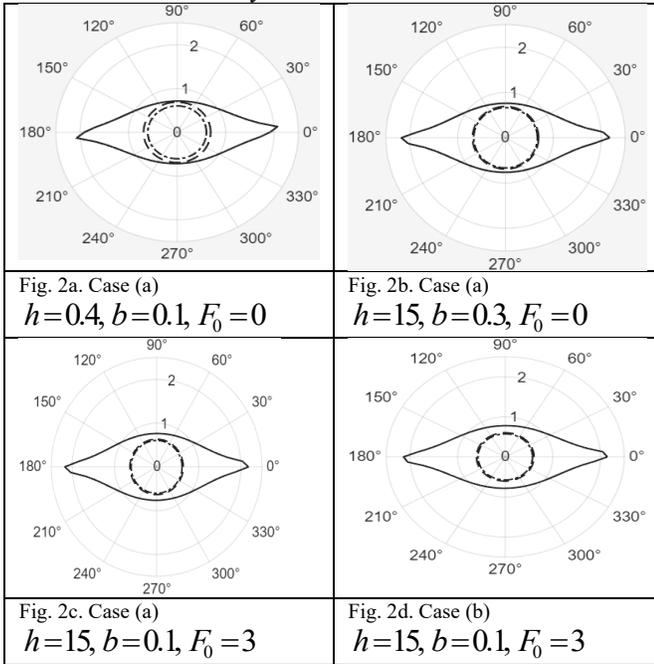

Fig. 2a. Case (a) $h=0.4, b=0.1, F_0=0$
Fig. 2b. Case (a) $h=15, b=0.3, F_0=0$
Fig. 2c. Case (a) $h=15, b=0.1, F_0=3$
Fig. 2d. Case (b) $h=15, b=0.1, F_0=3$

**Figure 2**. Simulated boundaries of the trapping ($F_0 > 0$) and stability ($F_0 = 0$) regions for equations (5.1)–(5.3).

Simulations of equations (5.2) and (5.3) efficiently estimate embedded disk radii that remain close to the minimal dimensions of the true regions and exhibit weak parameter sensitivity. The arrangement of the curves in Figure 2 is consistent with theoretical inferences.

## VI. Conclusion

Assessing boundedness, stability, and the temporal evolution of solution norms in time-varying vector nonlinear DDEs remains challenging, since the existing methods are often conservative and may overlook essential system dynamics.

We propose a rigorous yet computationally tractable framework based on constructing scalar auxiliary equations whose solutions upper-bound the norms of solutions of the original vector systems under aligned history functions. This reduction transforms the analysis of complex vector DDEs into study of scalar counterparts, enabling simplified analytical treatment and efficient, rigorous simulations.

This framework yields new boundedness, stability and robust stability criteria, along with estimates of the radii of the balls of history-function leading to bounded or stable solutions. Applying a Lipschitz-type condition to the nonlinear terms yields linear scalar equations that upper-bound their nonlinear analogs, thereby extending and simplifying stability and boundedness analysis, which, for example, brings new and less conservative ISS criterion.

A valuable advantage of the proposed method is that it naturally aggregates parameters, history functions, and uncertainty models through a unified framework, thereby simplifying subsequent sensitivity and robustness analyses.

Representative simulations show that solutions of the scalar auxiliary equation closely track the solution norms of the original vector systems initialized within the core of the boundedness and stability regions. The method also provides effective estimates of the radii of history-function balls that ensure bounded or stable behavior over a broad parameter range.

Limitations include the computational effort required to obtain the fundamental solution matrix of a linear nonautonomous system, the conservatism introduced by applying the norm inequalities and potential long-term growth of the scalar function $c(t)$. These issues will be addressed in future work under additional assumptions that preserve practical value of the results. Further research will also develop recursive approximations for the boundaries of trapping and stability regions as well as bilateral bounds on solution norms in vector nonlinear and nonautonomous delay systems.

**Appendix A**.
**Lemma 2**. Consider two scalar DDEs:
$$\dot{u}_i = f_i(t, u_i(t), u_i(t-h_1), ..., u_i(t-h_m)), \forall t \geq t_0;$$
$$u_i(t) = \varphi_i(t), \forall t \in [t_0 - \bar{h}, t_0], i=1,2$$
where $u_i(t,\varphi) \in \mathbb{R}, \forall t \geq t_0, \forall \|\varphi_i\| < \bar{\varphi} \in \mathbb{R}_+$ are unique solutions of the above DDEs, $h_i \in C([t_0,\infty); \mathbb{R}_+)$ $f_i \in C([t_0,\infty) \times \mathbb{R}^{m+1}; \mathbb{R})$, $\varphi_i \in C([t_0-\bar{h}, t_0]; \mathbb{R}), i=1,2$.

Assume also that $\varphi_1(t) \leq \varphi_2(t), \forall t \in [t_0-\bar{h}, t_0]$,
$$f_1(t, x_1, x_2, ..., x_{m+1}) \leq f_2(t, x_1, \bar{x}_2, ..., \bar{x}_{m+1}),$$
$\forall t \geq t_0, \forall x_j \leq \bar{x}_j \in \mathbb{R}, j=2,...,m+1$ and
$$\sup h_i(t) \leq \bar{h} < \infty, t \geq t_0; \inf h_i(t) \geq \underline{h} > 0, t \geq t_0, i=1,...,m \quad (A.1)$$
Then, $u_1(t,\varphi) \leq u_2(t,\phi), \forall t \geq t_0$.

**Proof**. Condition (A.1) allows time-stepping with a step size at least $h$, which reduces the problem on each successive interval to the corresponding well-known statement for scalar ODEs [32]. Consequently, the proof proceeds by induction. □

**Appendix B**. Assume that $x = [x_1 \ x_2]^T$ and
$$f = \begin{bmatrix} a_1(t)x_1^3 x_2^2(t-h_1) & a_2(t)x_2^3(t-h_2) \end{bmatrix}^T.$$
Then
$$|f|_2 \le |f|_1 \le |a_1||x_1^3||x_2^2(t-h_1)| + |a_2||x_2^3(t-h_2)|$$
$$\le |a_1||x(t)|^3|x(t-h_1)|^2 + |a_2||x(t-h_2)|^3$$
since $|x_i^n| \le |x|^n$, $|x_i^n(t-h)| \le |x(t-h)|^n$, $i=1,2$, $n \in \mathbb{N}$. Clearly, such reasoning can be extended on power series and some rational functions under the appropriate conditions.

**Acknowledgement**. The code for simulations discussed in Section V was developed by Steve Koblik.